\documentclass [11pt]{amsart}
\usepackage{bbm}
\usepackage{amssymb}
\usepackage{epic}
\usepackage{ecltree}
\usepackage{tikz}
\numberwithin{equation}{section}
\usepackage[normalem]{ulem}
\newtheorem{theorem}{Theorem}[section]
\newtheorem{lemma}[theorem]{Lemma}

\newtheorem{corollary}[theorem]{Corollary}

\newtheorem{question}[theorem]{Question}
\newtheorem*{theoremA}{Theorem A}
 \newtheorem*{theoremB}{Theorem B}
\newtheorem*{theoremC}{Theorem C}
\newtheorem*{theoremD}{Theorem D}
\theoremstyle{definition}
\newtheorem{definition}[theorem]{Definition}

\newtheorem{remark}[theorem]{Remark}
\newtheorem*{acknowledgement}{Acknowledgement}

\begin{document}
\title[$2$-rotund norms and $1$-unconditional and $1$-symmetric bases]{$2$-rotund norms for unconditional \\ and  symmetric sequence spaces}
%\dedicatory{}

%------------------------------------------------------------------------------------------------------------------------{Theorem C}------------------------

%------------------------------------------------------------------------------------------------------------------------------------------------

\author[Stephen Dilworth]{Stephen Dilworth}

\address{Department of Mathematics, University of South Carolina, Columbia,
SC 29208, USA}

\email{dilworth@math.sc.edu}

%------------------------------------------------------------------------------------------------------------------------------------------------

\author{Denka Kutzarova}

\address{Department of Mathematics, University of Illinois Urbana-Champaign,
Urbana, IL 61807, USA; Institute of Mathematics and Informatics, Bulgarian Academy of Sciences, Sofia, Bulgaria}

\email{denka@illinois.edu}

%------------------------------------------------------------------------------------------------------------------------------------------------
\author[Pavlos  Motakis]{Pavlos Motakis}

\address{Department of Mathematics and Statistics, York University, 4700 Keele Street, Toronto, Ontario, M3J 1P3, Canada}

\email{pmotakis@yorku.ca}

\thanks{ The first author was supported by Simons Foundation Collaboration Grant No. 849142.
 The second author was supported by Simons Foundation Collaboration Grant No. 636954.  The third author was supported by NSERC Grant RGPIN-2021-03639.}

%\date{\today}



\begin{abstract} A reflexive Banach space with an unconditional basis admits an equivalent  $1$-unconditional $2R$ norm and embeds into a reflexive space with a $1$-symmetric $2R$ norm. Partial results on $1$-symmetric  $2R$ renormings of spaces with a symmetric basis are obtained.
\end{abstract}

\maketitle  
\section{Introduction} \label{sec: Intro} 
The  notions of $2$-rotund and weakly $2$-rotund norms   were introduced by Milman \cite{M} and are defined as follows.
 \begin{definition}   Let $X$ be a Banach space. We say that a norm $\|\cdot\|$ on $X$ is $2$-rotund  ($2R$) (resp.\, weakly $2$-rotund ($W2R$))   if for every $(x_n) \subset X$ such that $\|x_n\| \le 1$ ($n \ge 1$) and
$$ \lim_{m,n \rightarrow \infty} \|x_m + x_n\| = 2,$$
there exists $x \in X$ such that $x = \lim_{n \rightarrow \infty} x_n$ strongly (resp.\, weakly). \end{definition}

Note that a $W2R$ norm is strictly convex.
It follows from a characterization of reflexivity due to  James \cite{J} that if  $X$ admits an equivalent $W2R$ norm then $X$ is reflexive. H\'ajek and Johanis 
\cite{HJ} proved the converse:   every reflexive Banach space admits an equivalent $W2R$ norm.  
 Odell and Schlumprecht \cite{OS} proved that every separable   reflexive Banach space $X$ admits an equivalent $2R$ norm (cf.\,\cite{G}). However,  it is an open question whether every  reflexive Banach space admits an equivalent 
$2R$ norm (cf. \cite[p. 72]{HJ}).

One  motivation for the present  article is the following result  of Figiel and Johnson  which combines Theorem 3.1 and Remark 3.2 of \cite{FJ}. 
Only the unconditional case  is stated explicitly in \cite{FJ},  but the argument for the unconditional case  also proves the symmetric case. 

\begin{theoremA} Let $X$ be a superreflexive Banach space with an unconditional (respectively, symmetric)  basis $(e_n)_{n=1}^\infty$. Then  $X$ admits an equivalent uniformly convex norm for which $(e_n)_{n=1}^\infty$ is $1$-unconditional (respectively, $1$-symmetric).
\end{theoremA}

Enflo \cite{E} showed that a space is superreflexive if and only if it admits an equivalent uniformly convex norm.  By the  theorem of Odell and Schlumprecht
above  a
separable space is reflexive  if and only if it admits an equivalent  $2R$ norm. Therefore it is natural to ask whether the analogue of Theorem~A holds for  $2R$ renormings of separable reflexive spaces. 

In Section~3 we prove  the analogous result in the unconditional case: a reflexive space with an unconditional basis  admits a $1$-unconditional  $2R$ norm. For the symmetric case, however, we have only partial results. In particular, the following question is open.

\begin{question}
\label{question2}
Let $X$ be a reflexive  Banach space with a symmetric basis $(e_n)_{n=1}^\infty$. Does $X$ admit an equivalent $2R$ norm for which 
$(e_n)_{n=1}^\infty$ is $1$-symmetric?
\end{question}

We show that the answer is positive if the lower Boyd index $p_X$ of $X$ satisfies $p_X>1$. We also prove that if $X$ is a reflexive space with an unconditional basis then $X$ is isomorphic to a $1$-complemented subspace of a space with a $2R$ norm and a $1$-symmetric basis. This  is a refinement of  a theorem of { Szankowski \cite{Sz}.
A similar argument proves that the non-superreflexive space with a symmetric basis which does not contain $c_0$ or $\ell_p$  constructed in \cite{FJ} and the space with a unique symmetric basic sequence (not equivalent to the unit vector basis of $c_0$ or $\ell_p$) constructed by Altshuler \cite{A} both  admit an equivalent $2R$ norm for which the basis is $1$-symmetric. 

A second motivation   is the open question  whether nonseparable reflexive spaces admit a $2R$ norm. To attack this problem  it is natural to examine specific classes of nonseparable reflexive spaces for  potential counterexamples or positive results.   In \cite{DK}  partial positive  results were obtained for nonseparable generalized Baernstein spaces. Another natural class  to examine is the  class of  spaces with an uncountable symmetric basis.
\begin{question}
\label{question3}
Suppose $X$ is a  reflexive  space with an uncountable symmetric basis $(e_\gamma)_{\gamma \in \Gamma}$. Does $X$ admit an equivalent
(not necessarily $1$-symmetric) $2R$ norm? \end{question} 

A positive answer to Quesion \ref{question2} would imply a positive answer to Question \ref{question3}.  In fact, it would imply that  $X$ admits a $1$-symmetric $2R$ norm.

In the final section  we show that $L_\infty[0,1]$ admits an equivalent rearrangement-invariant norm  which restricts to a $W2R$ norm on every reflexive subspace.

Finally, let us mention that related results are proved in \cite{GH} and \cite{HZ}. In \cite{GH} it is proved that a uniformly smooth (resp. uniformly convex)
space with a Schauder basis admits a uniformly smooth (resp. uniformly convex) renorming for which the basis is monotone, while in \cite{HZ}  the spaces with a symmetric basis which admit equivalent symmetric norms that are G\^{a}teaux differentiable or uniformly rotund in every direction are characterized.
\section{Preliminary results}
We shall use the following  characterization of $2$-rotundity (see e.g., \cite[II.6.4]{DGZ} or \cite{HJ}): $\|\cdot\|$ is a $2R$ norm on $X$  if for all  $(x_n)_{n=1}^\infty \subset X$ such that \begin{equation} \label{eq: alternativedef}
\lim_{m,n \rightarrow \infty} [\| x_m + x_n\|^2 - 2(\|x_m\|^2 + \|x_n\|^2)] = 0, \end{equation}
 there exists $x \in X$ such that $x=\lim_{n \rightarrow \infty} x_n$ strongly.

Day \cite{D}  introduced the  norm $\|\cdot \|_{\textrm{Day}}$ on $c_0$  defined by 
$$\|(a_n)_{n=1}^\infty\|_{\textrm{Day}} = (\sum_{n=1}^\infty 4^{-n} a_n^{*2})^{1/2},$$
where $(a_n^*)_{n=1}^\infty$ is the non-increasing rearrangement of $(|a_n|)_{n=1}^\infty$.
Let $(Y,\|\cdot\|)$ be a reflexive Banach space with normalized basis $(e_n)_{n=1}^\infty$. We define an equivalent norm
on   $Y$ thus:
\begin{equation}\label{eq: w2Rrenorming}   |\!|\!|\sum_{n=1}^\infty a_n e_n  \  |\!|\!| =( \|\sum_{n=1}^\infty a_ne_n\|^2 +  \|(a_n)_{n=1}^\infty\|_{\textrm{Day}}^2)^{1/2} .  
\end{equation}

 We will use the  following result of  H\'ajek and Johanis.  It is a consequence of Theorem~3 and Corollary~4 of \cite{HJ} and the reflexivity of $Y$. (Here $\|\sum_{n=1}^\infty a_n e_n\|_\infty  = \sup_{n \ge 1} |a_n|$ as usual.)
\begin{theoremB}  \label{lem: HJresult} 

Suppose $(y_n)_{n=1}^\infty \subset Y$ satisfies \begin{equation} \label{eq: triplenormcondition_2}
 \lim_{m,n \rightarrow \infty}[ |\!|\!|y_n + y_m |\!|\!|^2- 2(|\!|\!|y_n|\!|\!|^2 + |\!|\!|y_m|\!|\!|^2) ]= 0. \end{equation}
Then there exists $y \in Y$ such that \begin{equation*}
 y_n\rightarrow y\quad \text{weakly as $n \rightarrow \infty$}\end{equation*} and
$$ \lim_{n \rightarrow \infty} \|y_n - y\|_\infty = 0.$$
\end{theoremB}
For $K \ge 1$,  a basis $(e_n)_{n=1}^\infty$  is  $K$-unconditional if 
$$\|\sum_{n=1}^\infty \pm  a_n e_n\| \le K\|\sum_{n=1}^\infty a_n e_n\|$$ for all scalars $(a_n)_{n=1}^\infty$ and all choices of signs. The basis is $K$-symmetric if 
$$\|\sum_{n=1}^\infty \pm  a_{\sigma(n)} e_n\| \le K\|\sum_{n=1}^\infty a_n e_n\|$$
for all scalars $(a_n)_{n=1}^\infty$, all choices of signs, and all permutations $\sigma\colon \mathbb{N} \rightarrow \mathbb{N}$.

We refer the reader to \cite{LT1} for other  unexplained  Banach space  notation and terminology.
\section{$1$-unconditional bases}
\begin{theorem}\label{thm: 2Runconditional}  Suppose  that $X$  has an unconditional basis. Then $X$ admits an equivalent $1$-unconditional norm $ |\!|\!|\cdot|\!|\!|$ such that if $(x_n)_{n=1}^\infty\subset X$ is relatively weakly compact and satisifies \eqref{eq:  alternativedef}, then $(x_n)$ converges strongly.
In particular, if $X$ is reflexive, then $ |\!|\!|\cdot|\!|\!|$ is $2R$ and $1$-unconditional.
\end{theorem} 
\begin{proof} The proof closely  follows   \cite[Main Theorem]{OS}. So it  suffices to indicate how to adapt \cite[Main Theorem]{OS} to produce  a $1$-unconditional basis as well as a $2R$ norm. 

 Let $(e_n)_{n=1}^\infty$ be a semi-normalized unconditional basis for $X$ and let
 $\|\cdot\|$ denote any equivalent norm on $X$ which is strictly convex and for which $(e_n)_{n=1}^\infty$ is $1$-unconditional.
To see that such a norm exists, let 
$|\cdot|$ be any equivalent norm on $X$. Let
$$\|\sum_{n=1}^\infty a_n e_n\| := \sup |\sum_{n=1}^\infty \pm  a_n e_n| + (\sum_{n=1}^\infty 2^{-4n} a_n^2)^{1/2},$$
where the supremum is taken over all choices of signs. Then $\|\cdot\|$ is strictly convex and $(e_n)_{n=1}^\infty$ is a $1$-unconditional basis for $(X,\|\cdot\|)$.
For $x \in  X$, following \cite[p. 148]{OS}, define an equivalent norm $\|\cdot\|_x$ thus:
$$\|y\|_x := \|\|y\|x + y\| + \|\|y\|x - y\| \qquad (y \in X).$$
Then (\cite[Lemma 2.1]{OS}) \begin{equation} \label{eq: equivalentnorm} 2\|y\| \le \|y\|_x \le (2 + 2\|x\|)\|y\|.\end{equation}

Let $C$  be the countable vector space over $\mathbb{Q}$ defined by
$$C := \{ \sum a_n e_n \colon (a_n) \in c_{00}, a_n \in \mathbb{Q}, n \ge 1\}.$$
 Let us  say that $c = \sum a_ne_n \in C$ and $d = \sum b_ne_n \in C$  are absolutely equivalent if $|a_n| =|b_n|$ for all $n\ge 1$.
Note that absolute equivalence is an equivalence relation on $C$ and that the equivalence classes are finite. Let $\mathcal{A}$ be the collection of equivalence classes. For all  $A \in \mathcal{A}$ and for all  absolutely equivalent  $y,z \in C$,
 by $1$ unconditionality of $\|\cdot\|$ and absolute equivalence of $y$ and $z$, we have \begin{align} \begin{split} \label{align: uncond}
\sum_{c \in A} \|y\|_c  &= \sum_{c \in A}( \|\|y\|c + y\| + \|\|y\|c - y\|)\\
&= \sum_{c \in A}( \|\|z\|c + z\| + \|\|z\|c - z\|)\\
&=\sum_{c \in A} \|z\|_c.
\end{split}\end{align} Choose $p_A>0$ ($A \in \mathcal{A}$) such that  $\sum_{A \in \mathcal{A} }p_A(1 + \sum_{c \in A} \|c\|) < \infty$.
Define a norm $ |\!|\!|\cdot|\!|\!|$ on $X$ as follows:
 $$ |\!|\!|x|\!|\!| = \sum_{A \in \mathcal{A}} p_A \sum_{c \in A} \|x\|_c.
$$ It follows from \eqref{eq: equivalentnorm} that  $ |\!|\!|\cdot|\!|\!|$ is an equivalent norm.  Note  that  $ |\!|\!|\cdot|\!|\!|$
is strictly convex since  $\|\cdot\|_0 = 2\|\cdot\|$ which is strictly convex. Suppose $y,z \in C$ are absolutely equivalent. Then  \eqref{align: uncond} implies that
$$ |\!|\!|y|\!|\!| = \sum_{A \in \mathcal{A}} p_A \sum_{c \in A} \|y\|_c =  \sum_{A \in \mathcal{A}} p_A \sum_{c \in A} \|z\|_c =   |\!|\!|z|\!|\!|.$$
Since $C$ is dense in $X$, it follows that $(e_n)$ is a $1$-unconditional basis for $(X, |\!|\!|\cdot|\!|\!|)$.
The proof of \cite[Main Theorem]{OS},  especially Lemmas 2.2(a), 2.3(a), and 2.4, now shows that $\|\cdot\|_M$ is $2R$.
\end{proof}
\begin{remark} We are unable to adapt  the proof of Theorem~\ref{thm: 2Runconditional} to the case of a symmetric basis. The natural approach    would be to  say that  vectors from $C$  are  equivalent if their coefficient sequences are permuted. However, in this case the equivalence classes are infinite, so the proof does not go through. 
\end{remark}

Lindenstrauss \cite{L} proved that every space $X$ with an unconditional basis is isomorphic to a complemented subspace of a  space $Y$ with a symmetric basis.  Subsequently,   Szankowski \cite{Sz} proved that if $X$ is reflexive then  $Y$ can be chosen to be reflexive
and Davis \cite{D} proved that if $X$ is superreflexive then $Y$ can be chosen to be superreflexive. Davis's method also proves the reflexive case. As an application of Theorem~\ref{thm: 2Runconditional} we  use Davis's method to  prove the following  refinement  of Szankowski's  result.
\begin{theorem}\label{thm: Davisgeneralization}  Suppose that $X$ is reflexive and  has an unconditional basis. Then $X$ is isomorphic to a $1$-complemented subspace  of a space with a $1$-symmetric basis and a  $2R$ norm. Moreover, that subspace has a $1$-unconditional basis.
\end{theorem} We recall the presentation of  Davis's approach in \cite{LT1}[p. 125]. Let $(E, \|\cdot\|)$ and $(F,\|\cdot\|)$ be two Banach sequence spaces such that $(e_n)_{n=1}^\infty$ is a $1$-symmetric basis for both $E$ and $F$.  We assume also that $\|x\| _E \le \|x\|_F$ for all $x \in F$ and that
$$\lim_{n \rightarrow \infty} \frac{\|\sum_{i=1}^n e_i\|_E}{\|\sum_{i=1}^n e_i\|_F} = 0.$$ For each $m \ge 1$, define a $1$-symmetric  norm $\|\cdot\|_m$  on $E$
as follows:
$$ \|x\|_m = \inf \{ (\|y\|_E^2 + \|z\|_F^2)^{1/2} \colon x = \frac{1}{m} y + mz, y \in E, z \in F\}.
$$ Then (\cite[p. 125]{LT1}) \begin{equation}\label{eq: normestimates} \frac{1}{m} \|x\|_m \le \|x\|_E  \le 2m \|x\|_m \qquad(x \in E). \end{equation}
Hence $\|\cdot\|_m$ is equivalent to $\|\cdot\|_E$.

Now suppose that $X$ is a Banach space with a normalized $1$-unconditional basis $(f_n)_{n=1}^\infty$. For every strictly increasing sequence $(m_n)_{n=1}^\infty$ such that 
$\sum_{n=1}^\infty 1/m_n < \infty$, we define  the  space  $Y := Y(E,F,X, (m_n)_{n=1}^\infty)$  to be the collection of all $x \in E$ for which
\begin{equation}\label{eq: Ynormdef}  \|x\|_Y := \|\sum_{n=1}^\infty \|x\|_{m_n} f_n\|_X < \infty. \end{equation}  Then  $(e_n)$ is a $1$-symmetric basis for $Y$. (The condition $\sum_{n=1}^\infty 1/m_n < \infty$ guarantees that  $F$ embeds continuously  into $Y$ and, in particular,   that $(e_n)_{n=1}^\infty \subset Y$ \cite[p. 126]{LT1}.)

\begin{theoremC} \label{thm: davis} \cite[Prop. 3.b.4]{LT1}  For every $E,F$ and $X$ as above there exists an increasing sequence of numbers $(m_n)_{n=1}^\infty$ with $\sum_{n=1}^\infty 1/m_n< \infty$ such that $Y = Y(E,F,X,(m_n)_{n=1}^\infty)$  contains a complemented subspace isomorphic to $X$.
\end{theoremC} The following lemma  generalizes \cite[Lemma 3.b.11]{LT1}.
 \begin{lemma}\label{lem: subseq}  Suppose  that $E=c_0$ and $\sum_{n=1}^\infty1/ m_n < \infty$.  Let 
$$u_n = \sum_{i = q_n +1}^{q_{n+1}} c_i e_i \qquad(m \ge 1)$$  be a normalized block basis (with respect to $(e_n)_{n=1}^\infty$)  in $Y(c_0,F,X,(m_n)_{n=1}^\infty)$   such that $\lim_{i \rightarrow \infty} c_i = 0$.
Then some subsequence of $(u_n)_{n=1}^\infty$ is equivalent to a block basis of $(f_n)_{n=1}^\infty$ in $X$.
\end{lemma} \begin{proof}  Fix $N \ge 1$ and $n \ge 1$. It follows from \eqref{eq: normestimates} that 
$$\sum_{i=1}^N \|u_n\|_{m_i} \le (\sum_{i=1}^N m_i)\|u_n\|_E = (\sum_{i=1}^N m_i) \max \{|c_i| \colon q_n+1 \le i \le q_{n+1}\}. $$
Since $\lim_{i \rightarrow \infty} c_i = 0$,  we can inductively define an increasing sequence of natural numbers $1=N_1 < N_2<\cdots$ and a subsequence 
$(u_{n_k})_{k=1}^\infty$ such that for all $ k \ge 1$  $$\|\sum_{i=1}^{N_k} \|u_{n_{k}}\|_{m_i} f_i\|_X + \|\sum_{N_{k+1}+1}
^\infty \|u_{n_{k}}\|_{m_i} f_i\|_X < 2^{-k-1}.$$ It follows that  the block basis $(u_{n_k})_{k=1}^\infty \subset Y$ is equivalent to the block basis 
$$(\sum_{i = N_k+1}^{N_{k+1}} \|u_{n_k}\|_{m_i} f_i )_{k=1}^\infty \subset X.$$
\end{proof}
The next  lemma  is more general than is needed for the proof of Theorem~\ref{thm: Davisgeneralization}, but  we believe that the additional generality  may be  of independent interest.
 \begin{lemma}
  Suppose that $E=c_0$ and $X$ is reflexive. If $\sum_{n=1}^\infty 1/m_n < \infty$ then  $Y(c_0,F,X, (m_n)_{n=1}^\infty)$  is reflexive.
\end{lemma} \begin{proof}   Since $Y$ has a symmetric (hence unconditional) basis it follows from a result of James \cite{J2} that  $Y$ is reflexive unless $Y$ contains a subspace isomorphic to $c_0$ or $\ell_1$. 
We use the fact that every  (infinite-dimensional) subspace of $Y$ contains a further subspace  isomorphic to a subspace of $X$ or to a subspace of $E$ (see \cite[p.127]{LT1}).
 Since every subspace of $\ell_1$ contains a further subspace isomorphic to $\ell_1$,
and since neither $c_0$ nor $X$ contain a subspace isomorphic to $\ell_1$, it follows that $Y$ does not contain a subspace isomorphic to $\ell_1$. Suppose, to obtain a contradiction,  that $Y$ contains a sequence equivalent to the unit vector basis of $c_0$.  Since the unit vector basis of $c_0$ is weakly null, a standard gliding hump argument shows that $(e_n)_{n=1}^\infty$  admits a block basis 
$$u_n= \sum_{i = q_n +1}^{q_{n+1}} c_i e_i$$
equivalent to the unit vector basis of $c_0$. 
By the proof of \cite[Lemma 3.b.3]{LT1}, $$\text{$\max_{n\ge 1} \|\sum_{i=1}^n e_i\|_m  \rightarrow \infty$ as $m \rightarrow \infty$.}$$
So, using \eqref{eq: Ynormdef},  $\|\sum_{i=1}^n e_i\|_Y \rightarrow \infty$ as $n \rightarrow \infty$. 
 It   follows from the unconditionality  of $(e_n)_{n=1}^\infty$ and the fact that $\sup_{n\ge1}\|\sum_{i=1}^n u_i\|_Y < \infty$  that $\lim_{i \rightarrow \infty} c_i =0$. By Lemma~\ref{lem: subseq}, some subsequence of $(u_n)_{n=1}^\infty$
is equivalent to a block basis of $(f_n)_{n=1}^\infty \subset X$. Hence  $c_0$ is isomorphic to a subspace of $X$, which contradicts the reflexivity of $X$.
 \end{proof} \begin{remark}  Reflexivity of $Y(c_0, \ell_1,X,( 2^n)_{n=1}^\infty)$ was proved in \cite{FJ} using  results from \cite{DFJP}.
\end{remark}
  \begin{proof}[Proof of Theorem~\ref{thm: Davisgeneralization}] By Theorem~\ref{thm: 2Runconditional}, $X$ has an equivalent $2R$ norm $\|\cdot\|_X $  for which  $(f_n)_{n=1}^\infty$ is 
a $1$-unconditional basis. Let $Y:=Y(c_0, F, X, (m_n)_{n=1}^\infty)$ be as in Theorem~C.
We equip $Y$ with the equivalent norm defined by \eqref{eq: w2Rrenorming}. 

Suppose $(y_n)_{n=1}^\infty$ satisfies \eqref{eq: triplenormcondition_2}. Since $Y$ is reflexive,  by Theorem~B there exists $y \in Y$ such that  $y_n \rightarrow y$ weakly  and $\|y_n - y\|_\infty \rightarrow 0$ as $n \rightarrow \infty$.  Moreover,  \eqref{eq: triplenormcondition_2} implies that \begin{equation} \label{eq: Ycondition}
 \lim_{m,n \rightarrow \infty}[\|y_n + y_m \|_Y^2- 2(\|y_n\|_Y^2 + \|y_m\|_Y^2)] = 0. \end{equation}
Let $x_n = \sum_{i=1}^\infty \|y_n\|_{m_i} f_i$ ($n \ge1$).  By definition of $\|\cdot\|_Y$,  $\|x_n\|_X = \|y_n\|_Y$ . Moreover, by $1$-unconditionality of the basis  $(f_n)_{n=1}^\infty$ of 
$X$,  \begin{align*}
\|x_n + x_k\|_X &= \| \sum_{i=1}^\infty (\|y_n\|_{m_i} + \|y_k\|_{m_i}) f_i\|_X\\
& \ge \| \sum_{i=1}^\infty \|y_n\ + y_k\|_{m_i}\ f_i\|_X\\
&= \|y_n + y_k\|_Y.
\end{align*} Hence  \eqref{eq: Ycondition} implies that 
 \begin{equation*}
 \lim_{m,n \rightarrow \infty}[\|x_n + x_m \|_X^2- 2(\|x_n\|_X^2 + \|x_m\|_X^2)] = 0. \end{equation*}
Since $\|\cdot\|_X$ is $2R$, it follows that $(x_n)_{n=1}^\infty$ is a Cauchy sequence in $X$. Hence, given $\varepsilon >0$, there exists $N_1\ge 1$ such that
$$ \| \sum_{i=N_1+1}^\infty  \|y_n\|_{m_i} f_i\|_X < \frac{\varepsilon}{4} \qquad (n \ge 1).$$  Hence
 \begin{equation} \label{eq: Cauchy1}  \| \sum_{i=N_1+1}^\infty  \|y_n-y_k\|_{m_i} f_i\|_X < \frac{\varepsilon}{2} \qquad (n,k \ge 1). \end{equation}
Recall that, by \eqref{eq: normestimates},  $\|\cdot\|_m$ is equivalent to $\|\cdot\|_\infty$ for all $m \ge 1$. Since $\lim_{n \rightarrow \infty} \|y_n - y\|_\infty= 0$,
it follows that $(y_n)_{n=1}^\infty$ is a Cauchy sequence in $\|\cdot\|_m$ for all $m \ge 1$. Hence there exists $N_2 \ge 1$ such that
\begin{equation} \label{eq: Cauchy2}  \|\sum_{i=1}^{N_1} \|y_n - y_k\|_{m_i} f_i\|_X < \frac{\varepsilon}{2} \qquad (n,k \ge N_2). \end{equation}
Combining \eqref{eq: Cauchy1} and \eqref{eq: Cauchy2}, we have $\|y_n - y_k\|_Y < \varepsilon$ for all $n,k \ge N_2$. So $(y_n)_{n=1}^\infty$ is a Cauchy sequence in $Y$ and hence $\lim_{n \rightarrow \infty} \|y_n - y\|_Y = 0$.

The proof of Theorem~C shows that $X$ is isomorphic to the closed linear span  $Z$ of disjointly supported constant coefficient vectors in $Y$.
Hence $Z$ has a $1$-unconditional basis and is the range of an averaging projection on $Y$. So $Z$ is $1$-complemented in $Y$. 
\end{proof} Let $T$ be the space introduced in \cite{FJ}
(the dual of the space  that does not contain $c_0$ or $\ell_p$ constructed by Tsirelson \cite{T}). It was proved in  \cite{FJ} that $Y(c_0,\ell_1,T, (2^n)_{n=1}^\infty)$ does not contain a subspace isomorphic to  $c_0$ or $\ell_p$.

 Let $d_{w,1}$ be the Lorentz sequence space corresponding to the weight sequence $w = (1/n)$. The  norm in $d_{w,1}$  is given by
$$\|\sum_{n=1}^\infty a_n e_n\| = \sum_{n=1}^\infty\frac{ a_n^*}{n}.$$ It was proved by Altshuler \cite{A} that $Y(c_0, d_{w,1},T,(2^n)_{n=1}^\infty)$ has a unique symmetric basic sequence which, moreover, is not equivalent to the unit vector basis of $c_0$ or $\ell_p$.

The proof of Theorem~\ref{thm: Davisgeneralization} also establishes the following result.  
 \begin{theorem} The spaces $Y(c_0,\ell_1,T, (2^n)_{n=1}^\infty)$ of \cite{FJ} and  $Y(c_0, d_{w,1},T,(2^n)_{n=1}^\infty)$ of \cite{A} both  have equivalent $2R$ norms with a $1$-symmetric basis.
\end{theorem}
\section {$1$-symmetric bases}
In this section $(X, \|\cdot\|)$ denotes a reflexive  Banach space with a symmetric basis $(e_n)_{n=1}^\infty$. 

Let us recall the definition of the lower Boyd index (\cite{B})  $p_X$ of $X$  (cf. \cite[p. 130]{LT2}).
For $m \in \mathbb{N}$, the linear operator  $D_m \colon X \rightarrow X$ is defined by 
$$D_m(\sum_{n=1}^\infty a_n e_n) = \sum_{n=1}^\infty a_n (\sum_{j=(n-1)m+1}^{nm} e_j).$$
The lower Boyd index $p_X$ is defined by
\begin{equation} \label{eq: Boyd} p_X = \sup_{m \ge 2}  \frac{\log m}{\log \|D_m\|} = \lim_{m \rightarrow \infty} \frac{\log m}{\log \|D_m\|}. \end{equation}

The following main result of this section is an immediate consequence of Theorem~\ref{thm: 2Rrenorming} proved below.
\begin{theorem}\label{thm: lowerBoyd} Suppose that $X$ is a reflexive Banach space with a symmetric basis such that  $p_X >1$. Then $X$ admits a $1$-symmetric $2R$ norm.
\end{theorem}
\cite[Prop. 2.b.7]{LT2}, which  characterizes when $p_X >1$, yields a geometrical formulation of Theorem~\ref{thm: lowerBoyd}.
\begin{corollary} Suppose that  $X$  is a reflexive Banach space with a symmetric basis which does not admit uniformly isomorphic copies of $\ell_1^n$ spanned by disjointly supported vectors with the same distribution.
Then $X$ admits a $1$-symmetric $2R$ norm.
\end{corollary}
 For $x  = \sum_{i=1}^\infty x(i)e_i \in X$, define the formal series
$$ \widehat{x}: = \sum_{n=1}^\infty (\frac{1}{n} \sum_{i=1}^n x^*(i)) e_n.$$ 

We prove the following lemma for the sake of completeness. More general results in the setting of rearrangement-invariant function spaces rather than
symmetric  sequence spaces are proved in \cite[Theorem 5.15]{BS}.
\begin{lemma} Suppose that $p_X >1$. Then there exists a constant $c > 0$ such that
$$\|\widehat{x}\| \le c \|x\| \qquad (x \in X).$$
\begin{proof} We may assume that $(e_n)_{n=1}^\infty$ is a $1$-symmetric basis of $X$. Let $1 < p < p_X$. It follows from \eqref{eq: Boyd} that there exists $A>0$ such that $\|D_m\| \le A m^{1/p}$ for all $m \ge 1$.
Consider $x = \sum_{n=1}^\infty x(n) e_n \in X$, where $(x(n))_{n=1}^\infty$ is a nonnegative decreasing sequence.
Define $f \colon (0,\infty) \rightarrow (0,\infty)$ by $f(t) = x(n)$ for $n\ge1$ and  $n-1 < t \le n$. Then
$$ \widehat{x}  = \sum_{n=1}^\infty (\int_0^1 f(tn) \, dt) e_n = \int_0^1 (\sum_{n=1}^\infty f(tn)e_n) \, dt$$ Hence \begin{align*}
\|\widehat{x}\| &\le  \int_0^1 \|\sum_{n=1}^\infty f(tn)e_n\| \, dt\\
&\le  \sum_{m=1}^\infty  2^{-m}  \|\sum_{n=1}^\infty f(2^{-m}n)e_n\| \\
\intertext{(by $1$-unconditionality)}
&= \sum_{m=1}^\infty 2^{-m} \|D_{2^m}(x)\|\\
&\le \sum_{m=1}^\infty 2^{-m}  (A 2^{m/p}) \|x\|\\
&= \frac{A}{ 2^{1-1/p }-1}\|x\|.
\end{align*}
\end{proof}
\end{lemma} Henceforth, we suppose that $p_X>1$,  and, using Theorem~\ref{thm: 2Runconditional},  that  $\|\cdot\|$ is $2R$, and that
$(e_n)_{n=1}^\infty$ is a symmetric $1$-unconditional basis for $\|\cdot\|$. Suppose that  $(e_n)_{n=1}^\infty$ is $K$-symmetric for $\|\cdot\|$.
Define a quasi-norm $ |\!|\!|\cdot |\!|\!|$ as follows: \begin{equation} \label{eq: defoftriplenorm}
 |\!|\!|x |\!|\!| = (\|\widehat{x}\|^2 + \|(x(n))_{n=1}^\infty\|_{\text{Day}}^2)^{1/2} \qquad (x = \sum_{n=1}^\infty x(n) e_n \in X).
\end{equation}
\begin{lemma} 
\label{lemma triple-bar-norm symmetry}
$ |\!|\!|\cdot |\!|\!|$  is a $1$-symmetric equivalent norm on $X$.
\end{lemma} \begin{proof} Clearly, $|\!|\!|\cdot |\!|\!|$ is a $1$-symmetric quasi-norm since $\widehat{x}$ depends only on $(x^*(n))_{n=1}^\infty$ and
$\|\cdot \|_{\text{Day}}$ is $1$-symmetric.
 For $x \in X$,
\begin{align*}  \frac{1}{K} \|x\| &= \frac{1}{K} \|\sum_{n=1}^\infty x(n) e_n\|\\
& \le  \| \sum_{n=1}^\infty x^*(n) e_n\| \\ \intertext{(since $(e_n)_{n=1}^\infty$ is a $K$-symmetric basis)}
&\le \|\hat{x}\|\\
\intertext{(since $x^*(n) \le \frac{1}{n} \sum_{i=1}^n x^*(i)$ and $(e_n)_{n=1}^\infty$ is a $1$-unconditional basis)}
& \le c \|x\|.
\end{align*} Since $\|\cdot\|_{\text{Day}}$ is equivalent to $\|\cdot\|_\infty$, it follows that $\|\cdot\|$ and  $|\!|\!|\cdot |\!|\!|$ are equivalent
quasi-norms. It remains to prove that  $|\!|\!|\cdot |\!|\!|$ is in fact a norm, i.e., that  $|\!|\!|\cdot |\!|\!|$ satisfies the triangle inequality.
It suffices to show that $x \mapsto \|\widehat{x}\|$ satisfies the triangle inequality.
Let $x,y \in X$. Note that 
\begin{equation}  \label{eq: triangle} \widehat{(x+y)}(n) \le \widehat{x}(n)+ \widehat{y}(n) \qquad (n \in \mathbb{N}). \end{equation} 
Since $(e_n)_{n=1}^\infty$ is a $1$-unconditional basis, it follows that  \begin{equation} \label{eq: triangleineq}
\|\widehat{x+y}\| \le \|\widehat{x} + \widehat{y}\| \le  \|\widehat{x}\| +\| \widehat{y}\|.\end{equation}
Hence  $x \mapsto \|\widehat{x}\|$  and $\|\cdot\|$ are equivalent norms on $X$. 

\end{proof} For $x \in X$ and $N,M \in \mathbb{N}$, define
$$ x\cdot 1_{[N,M]} := \sum_{n=N}^M x(n)e_n.$$ 
\begin{lemma} \label{lem: xy} For $x,y \in X$ and $N \in \mathbb{N}$,
$$ | \|\widehat{x}\cdot 1_{[N,\infty)}\| - \|\widehat{y}\cdot 1_{[N, \infty)}\|| \le c \|x-y\|.$$
\end{lemma}  \begin{proof} \eqref{eq: triangle}  yields
$$\|\widehat{x}\cdot 1_{[N,\infty)}\| \le \|\widehat{y}\cdot 1_{[N,\infty)}\| +  \|(\widehat{x-y})\cdot 1_{[N,\infty)}\|.$$
Hence
$$ \|\widehat{x}\cdot 1_{[N,\infty)}\| - \|\widehat{y}\cdot 1_{[N,\infty)}\| \le \|(\widehat{x-y})\cdot 1_{[N,\infty)}\| \le \|(\widehat{x-y})\| \le c\|x-y\| .$$
Interchanging $x$ and $y$ gives the result.
\end{proof}\begin{lemma}\label{lem: normafterN} Suppose that $x \in X$, $y_n \in X$ ($n \ge 1$),   $\|y_n\| \ge \delta>0$, $\lim_{n \rightarrow \infty} \|y_n\|_\infty=0$
 and $\min(\operatorname{supp} (y_n)) \rightarrow \infty$.
Then, for all $N \ge 1$, 
$$ \liminf_{n \rightarrow \infty} \|(\widehat{x+y_n})\cdot 1_{[N, \infty)} \| \ge \frac{\delta}{K}.$$
\end{lemma} \begin{proof} Let $\alpha>0$. Choose $N_1 \in \mathbb{N}$ so that $x^{\prime} := x 1_{[1,N_1]}$ satisfies $\|x - x^{\prime}\| < \alpha$.
Then, for all sufficiently large $n$, we have
$$ \min ( \operatorname{supp}(y_n)) > N_1 \ge \max (\operatorname{supp}(x^{\prime}).$$
 Hence, for all sufficiently large $n$ and for all  $i \ge 1$,
$$(\widehat{x^{\prime} + y_n})(i) \ge y_n^*(i).$$ So  for all $N \ge 1$ and for all sufficiently large $n$,\begin{align*}
\|(\widehat{x^\prime + y_n)} 1_{[N,\infty)}\| &\ge \|\sum_{i = N+1}^\infty y_n^*(i) e_i\| \\
&\ge \|  \sum_{i=1}^\infty y_n^*(i) e_i\| - N\|y_n\|_\infty\\
& \ge \frac{1}{K} \|y_n\| -N \|y_n\|_\infty \\& \ge \frac{\delta}{K} - N\|y_n\|_\infty.
\end{align*} Hence,  by  Lemma~\ref{lem: xy},   for all $N \ge 1$ and for all sufficiently large $n$,  \begin{align*}
\| (\widehat{x + y_n})\cdot 1_{[N,\infty)}\| &\ge \| (\widehat{x ^\prime+ y_n})\cdot 1_{[N,\infty)}\|- c\|x - x^\prime\|\\ &\ge \frac{\delta}{K} - N\|y_n\|_\infty -c\alpha.
\end{align*} Since  $\lim_{n \rightarrow \infty} \|y_n\|_\infty=0$ and  $\alpha>0$ is arbitrary  the result follows.
\end{proof} \begin{theorem} \label{thm: 2Rrenorming} $ |\!|\!|\cdot|\!|\!|$ is a $1$-symmetric $2R$ equivalent norm on $X$.
\end{theorem}

\begin{proof}
 Let us summarize the relevant progress we have made so far in this section. We used Theorem \ref{thm: 2Runconditional} to equip $X$ with an equivalent 2R norm $\|\cdot\|$ that is 1-unconditional but not necessarily 1-symmetric (see the paragraph before Lemma \ref{lemma triple-bar-norm symmetry}). We then defined the equivalent norm $ |\!|\!|\cdot|\!|\!|$, which we have shown in Lemma \ref{lemma triple-bar-norm symmetry} to be 1-symmetric. It remains to prove that it is 2R. Suppose $(x_n)_{n=1}^\infty \subset X$ satisfies \begin{equation} \label{eq: triplenormcondition}
 \lim_{m,n \rightarrow \infty}[ |\!|\!|x_n + x_m |\!|\!|^2- 2(|\!|\!|x_n|\!|\!|^2 + |\!|\!|x_m|\!|\!|^2) ]= 0. \end{equation}
By Theorem~B,  $(x_n)_{n=1}^\infty$ converges weakly to some $x \in X$ and $\lim_{n \rightarrow \infty} \|x_n - x\|_\infty = 0$.
Let $x_n = x + y_n$ and suppose that  $(y_n)_{n=1}^\infty$ does not converge to zero  in norm. Since $\lim_{n \rightarrow \infty}\|y_n\|_\infty = 0$, a  gliding hump and an  approximation  argument show, after passing to a subsequence and relabelling, that without loss of generality each $y_n$ has finite support, that 
$(y_n)_{n=1}^\infty $ is a block basis with respect to $(e_n)_{n=1}^\infty$, and hence $\min( \operatorname{supp}(y_n)) \rightarrow \infty$ as $n \rightarrow \infty$, and that $\|y_n\| > \delta >0$ ($n \ge 1$). 

It follows from \eqref{eq: triplenormcondition}  and the definition of $|\!|\!|\cdot  |\!|\!|$ in \eqref{eq: defoftriplenorm}   that \begin{equation*}
 \lim_{m,n \rightarrow \infty} \|\widehat{x_n + x_m} \|^2- 2(\|\widehat{ x_n} \|^2 + \|\widehat{x_m}\|^2) = 0. \end{equation*}
Note that $ \|(\widehat{x_n + x_m}) \|\le  \|\widehat{x_n} +\widehat{ x_m} \| $ by \eqref{eq: triangleineq}. Hence \begin{equation*}
 \lim_{m,n \rightarrow \infty} \|\widehat{x_n} + \widehat{x_m} \|^2- 2(\|\widehat{ x_n} \|^2 + \|\widehat{x_m}\|^2) = 0. \end{equation*}
Since $\|\cdot\|$ is a $2R$ equivalent norm on $X$, it follows that $(\widehat{x_n})_{n=1}^\infty$ converges strongly in $X$. 
By Lemma~\ref{lem:  normafterN}, for all $N \ge 1$,
$$\liminf_{m \rightarrow \infty} \|\widehat{x_m}\cdot 1_{[N,\infty)} \|= \liminf_{m \rightarrow \infty}\|(\widehat{x+y_m})\cdot 1_{[N,\infty)}\|  \ge \frac{\delta}{K},$$ 
which  contradicts the fact that $(\widehat{x_n})_{n=1}^\infty$ is a Cauchy sequence in $X$.
\end{proof} \section{Symmetric renormings of $\ell_\infty$ and $L_\infty$.}
 A symmetric  renorming of $\ell_\infty$ \cite{HZ} is an equivalent norm $\|\cdot\|$ on $\ell_\infty$  such that
$$ \|(a_n)_{n=1}^\infty\| = \|(a_{\sigma(n)})_{n=1}^\infty\| \qquad((a_n)_{n=1}^\infty \in \ell_\infty)$$
for  all permutations $\sigma$ of $\mathbb{N}$.
It was proved in \cite[Proposition~5]{HZ} that for  every symmetric renorming  $\|\cdot\|$, $(\ell_\infty,\|\cdot\|)$ contains a subspace isometric to $(\ell_\infty, \|\cdot\|_\infty)$.  For the sake of completeness we include an elementary proof that avoids uncountable cardinals.

\begin{theorem} \label{thm: symmetricrenormingoflinfty} Let  $\|\cdot\|$ be a $1$-symmetric norm on $\ell_\infty$. Then there exists a subspace $Y$ of $(\ell_\infty,\|\cdot\|)$  that is isometrically isomorphic to $(\ell_\infty,\|\cdot\|_\infty)$. In fact, there exists $\alpha>0$ such that,  for all $y \in Y$,  $\|y\| = \alpha \|y\|_\infty$.
\end{theorem}

\begin{proof}
Let  $\|\cdot\|$ be a $1$-symmetric norm on $\ell_\infty$. We let $2\mathbb{N}$ denote the set of even positive integers and $\ell_\infty(2\mathbb{N})$ the subspace of $\ell_\infty$ comprising all $x$ with $\mathrm{supp}(x) = \{i\in\mathbb{N}:x(i)\neq 0\}\subset2\mathbb{N}$. We will first show that $\|\cdot\|$ restricted on $\ell_\infty(2\mathbb{N})$ is 1-suppression unconditional, i.e., for every $x,y\in \ell_\infty(2\mathbb{N})$ such that, for all $i\in\mathbb{N}$, $x(i) = y(i)$ or $y(i) = 0$, we have $\|x\| \geq \|y\|$. We verify this on the dense linear subspace of $\ell_\infty(2\mathbb N)$ consisting of all $x$ that have the form
\[x = \sum_{i=1}^na_i\chi_{A_i},\]
where $n\in\mathbb N$, $a_1,\ldots,a_n$ are (not necessarily different) scalars, and $A_1,\ldots,A_n$ are disjoint (finite or infinite) subsets of $2\mathbb N$. By symmetry, it suffices to show that, letting
\[y = \sum_{i=1}^{n-1}a_i\chi_{A_i},\]
$\|x\| \geq \|y\|$. Fix an infinite $S\subset \mathbb N\setminus 2\mathbb N$. For each $N\in\mathbb N$, choose disjoint subsets $A_n^{1},\ldots,A_n^{N}$ of $S$ that are equinumerous to $A_n$, and, for $1\leq j\leq N$, let $x_N^{j} = y+a_n\chi_{A_n^j}$. By symmetry, $\|x_N^j\| = \|x\|$, and thus letting
\[y_N = \frac{1}{N}\sum_{j=1}^Nx_N^{j}\]
we have $\|x\| \geq \|y_N\|$. At the same time, $\|y_N - y\|_\infty\to 0$, and thus by equivalence $\|y_N - y\|\to 0$, which yields $\|x\| \geq \|y\|$.

By scaling and symmetry, we may assume that, for all $n \in \mathbb{N}$,  $\|e_n\|=1$ and thus, by 1-suppression unconditionality, for all $x \in \ell_\infty(2\mathbb N)$, $\|x\| \ge \|x\|_\infty$. Put
$$\alpha = \sup \{\|x\| \colon  x \in \ell_\infty(2\mathbb{N}),\|x\|_\infty = 1\}.$$
Then there exists $x_0 \in \ell_\infty(2\mathbb{N})$ such that $\|x_0\|_\infty = 1$ and $\|x_0\|= \alpha$. Indeed, by symmetry we can pick disjointly supported vectors $(x_n)_{n=1}^\infty$ in $\ell_\infty(2\mathbb{N})$ such that, for all $n  \in \mathbb{N}$, $\|x_n\|_\infty = 1$ and  $\|x_n\| \ge \alpha -  1/n$.
This follows from the symmetry of $\|\cdot\|$ and the fact that any bijection between infinite subsets of $2 \mathbb{N}$ extends to a permutation of $\mathbb{N}$. 
Define $x_0 = \sum_{n=1}^\infty x_n$ pointwise. Then, $\|x_0\|_\infty = 1$ and thus $\|x_0\| \le \alpha$. By 1-suppression unconditionality,  for all $n \in \mathbb{N}$, $\|x_0\| \geq \|x_n\| \geq \alpha-1/n$. So $\|x_0\| = \alpha$. Pick a disjointly supported sequence $(y_n)_{n=1}^\infty$ in $\ell_\infty(2\mathbb{N})$ so that each $y_n$ has the same distribution as $x_0$. Then, for all $(a_n) \in \ell_\infty$,     $\sum_{n=1}^\infty a_n y_n$ (defined pointwise) satisfies $\|\sum_{n=1}^\infty a_n y_n\| = \alpha \|(a_n)\|_\infty.$ Hence $Y = \{\sum_{n=1}^\infty a_n y_n \colon (a_n) \in \ell_\infty\}$ has the required property.
\end{proof}

On the other hand,  by \cite[Corollary 4]{HJ} $\ell_\infty$ admits an equivalent norm which restricts to a $W2R$ norm on reflexive subspaces. Clearly,  every such norm is strictly convex and hence cannot be symmetric by Theorem~\ref{thm: symmetricrenormingoflinfty}.

Next we consider   rearrangement-invariant renormings of $L_\infty[0,1]$. Curiously,  we  reach a  rather different conclusion  from the  case of  $\ell_\infty$.

 We will apply the following  result from \cite{DKT}.
\begin{theoremD} \cite{DKT} There is an equivalent rearrangement-invariant  (Orlicz) norm  $|\!|\!|\cdot  |\!|\!|$ on $L_1[0,1]$ satisfying the following  restricted uniform convexity  
condition. 
Let $K \subset \{x \in L_1[0,1] \colon  |\!|\!|x |\!|\!| \le 1\}$ be weakly compact. Then, given $\varepsilon>0$, there exists $\delta>0$ such that
for all $x,y \in K$, 
$$|\!|\!|x + y |\!|\!| > 2 - \delta \Rightarrow   |\!|\!|x - y |\!|\!| < \varepsilon.$$
\end{theoremD} 
\begin{corollary} Suppose  $(y_n)_{n=1}^\infty$ is relatively weakly compact in $L_1[0,1]$ and satisfies 
\begin{equation} \label{eq: triplenormcondition_3} 
 \lim_{m,n \rightarrow \infty}[ |\!|\!|y_n + y_m |\!|\!|^2- 2(|\!|\!|y_n|\!|\!|^2 + |\!|\!|y_m|\!|\!|^2)] = 0. \end{equation}
 Then
$(y_n)_{n=1}^\infty$ converges in $L_1[0,1]$.
\end{corollary} \begin{proof}  The proof is  omitted as it  is essentially  the same as the proof that a uniformly convex norm is $2R$.
\end{proof}
\begin{theorem}
 Let $(Y,|\cdot|)$ be a rearrangement-invariant space on $[0,1]$. Then $Y$ admits an equivalent rearrangement-invariant norm 
$\|\cdot\|$  such that if $(y_n)_{n=1}^\infty$ is relatively weakly compact in $Y$ and satisfies
\begin{equation} \label{eq: YW2R}  \lim_{m,n \rightarrow \infty}[ \|y_m + y_n\|^2 - 2(\|y_m\|^2 + \|y_n\|^2) ]= 0, \end{equation}
then  $(y_n)_{n-1}^\infty$ converges weakly in $Y$.
In particular,  $\|\cdot\|$ restricts to a  $W2R$ norm on every reflexive subspace of $Y$.
\end{theorem} \begin{proof}  Note that $Y$ embeds continuously into $L_1[0,1]$.  Define $\|\cdot \|$ as follows:
$$\|y\| = (|y|^2 + |\!|\!|y|\!|\!|^2)^{1/2} \qquad (y \in Y).$$ 
Suppose that $(y_n)_{n=1}^\infty$ satisfies \eqref{eq: YW2R}. Then $(y_n)_{n=1}^\infty$ also satisfies \eqref{eq: triplenormcondition_3} 
and $(y_n)_{n=1}^\infty$ is relatively weakly compact in $L_1[0,1]$.
It follows from Theorem~D that $(y_n)_{n=1}^\infty$ converges in $L_1[0,1]$, which implies that $(y_n)_{n=1}^\infty$ has a unique weak cluster point in $Y$, i.e. that
 $(y_n)_{n=1}^\infty$ converges weakly in $Y$. \end{proof}
\begin{corollary} $L_\infty[0,1]$ admits an equivalent rearrangement-invariant norm which restricts to a $W2R$ norm on every reflexive subspace.
\end{corollary}
 \begin{acknowledgement}  We thank the referees for bringing several references to our attention and for their insightful comments which improved  the presentation of the article.
\end{acknowledgement}

\end{document}